\documentclass[pdflatex,sn-mathphys]{sn-jnl}% Math and Physical Sciences Reference Style
\usepackage{color}
\jyear{2023}%

\theoremstyle{thmstyleone}%
\theoremstyle{thmstyletwo}%

\theoremstyle{thmstylethree}%

\begin{document}

\title[\textsc{LOBPCG}]{A robust, open-source implementation of the locally optimal block preconditioned conjugate gradient for large eigenvalue problems in quantum chemistry}


\author[1]{\fnm{Tommaso} \sur{Nottoli}}
\equalcont{These authors contributed equally to this work.}
\author[1]{\fnm{Ivan} \sur{Giann\`i}}
\equalcont{These authors contributed equally to this work.}
\author[2]{\fnm{Antoine} \sur{Levitt}}
\author[1]{\fnm{Filippo} \sur{Lipparini}*}\email{filippo.lipparini@unipi.it}
\affil[1]{\orgdiv{Dipartimento di Chimica e Chimica Industriale}, \orgname{Universit\`a di Pisa}, \orgaddress{\street{Via G. Moruzzi 13}, \city{Pisa}, \postcode{56124}, \country{Italy}}}

\affil[2]{\orgdiv{Laboratoire de Math{\'e}matiques d'Orsay}, \orgname{Universit\'e Paris-Saclay}, \orgaddress{\city{Orsay}, \postcode{91405}, \country{France}}}

% \affil[3]{\orgdiv{Department}, \orgname{Organization}, \orgaddress{\street{Street}, \city{City}, \postcode{610101}, \state{State}, \country{Country}}}


\abstract{We present two open-source implementations of the Locally
  Optimal Block Preconditioned Conjugate Gradient (\textsc{LOBPCG})
  algorithm to find a few eigenvalues and eigenvectors of large,
  possibly sparse matrices. We then test \textsc{LOBPCG} for various
  quantum chemistry problems, encompassing medium to large, dense to
  sparse, well-behaved to ill-conditioned ones, where the standard
  method typically used is Davidson's diagonalization. Numerical tests
  show that, while Davidson's method remains the best choice for most
  applications in quantum chemistry, LOBPCG represents a competitive
  alternative, especially when memory is an issue, and can even
  outperform Davidson for ill-conditioned, non diagonally dominant
  problems.}% \textbf{AL: water down ?}



\maketitle
\section{Introduction}\label{sec1}
Computing eigenvalues and eigenvectors of matrices is probably the
most prominent linear-algebra operation in quantum chemistry.
Eigenvalues and eigenvectors are computed in the self-consistent field
(SCF) algorithm,\cite{roothaan1,roothaan2,almloef_scf} in
restricted-step second-order optimizations,\cite{jensen1984direct}
response theory
calculations,\cite{christiansen1998response,casida1995time}, algebraic diagrammatic construction (ADC)\cite{schirmer1982beyond,dreuw2015algebraic}, and unitary coupled-cluster (UCC)\cite{bartlett1989alternative,taube2006new,liu2018unitary}. Moreover they are calculated in ground and
excited state configuration interaction (CI)
calculations,\cite{handy1980multi,olsen1988determinant} including in
state-average complete active space self-consistent field (CASSCF) \cite{werner2009matrix,roos1987complete,shepard1987multiconfiguration}. Quantum chemistry calculations are typically performed using localized basis sets, most commonly made by Gaussian atomic orbitals (GTOs). Thanks to the compactness of such basis sets, the generalized eigenvalue problem in SCF can be solved using dense linear-algebra techniques. Despite the cubic scaling, diagonalizing the Fock (or Kohn-Sham) matrix is usually not a significant bottleneck for standard calculations, which are instead dominated by the cost of assembling the Fock matrix. Iterative diagonalization techniques are therefore mostly used for post-Hartree Fock calculations, such as CASSCF and CI, where the combinatorially-scaling size of the CI Hamiltonian makes it impossible to build it in memory, let alone diagonalize it, but for the smallest cases. They are also used to compute the direction and step length in second-order SCF and CASSCF strategies, and to compute excited states in CI singles and time-dependent SCF.
Most of these applications involve computing one, or a small number
(up to a few hundreds) of eigenvalues and eigenvectors, and involve
symmetric, diagonally dominant matrices.

Since its introduction in 1975, the iterative method proposed by
Ernest R. Davidson has been the method of choice\cite{davidson1975}
and the default strategy used by the majority of quantum chemistry
codes. Davidson's diagonalization performs particularly well for
diagonally dominant symmetric eigenvalue problems, for which it
exhibits fast and robust convergence, and has been generalized to
non-symmetric problems and to deal with multiple eigenvalues and
eigenvectors \cite{liu1978simultaneous,zuev2015new}. 
However, the method suffers of a few drawbacks. First, its performance
degrades if the matrix is not diagonally
dominant. Even though most algorithms in quantum chemistry
involve diagonally dominant matrices, this is not always the case,
with the parameter Hessian in second-order CASSCF being a prominent
example. Finally, as the vector subspace used in
Davidson's method to expand the sought eigenvectors can become quite
large, especially if many eigenvalues are needed, the cost of
diagonalizing the matrix projected in the subspace and of
orthogonalizing the new trial vectors can become non-negligible, and
memory requirements problematic.

The Locally Optimal Block Preconditioned Conjugate Gradient (LOBPCG)
method offers a way out by more aggressively truncating the
approximation space, while attempting to preserve the convergence
properties of Davidson \cite{knyazev2001toward}. Since its
introduction, it has been successfully used in several applications,
and in particular in approaches to the quantum many-body problem in
condensed-matter physics \cite{bottin2008large} and
\cite{shao2018accelerating}, but to our knowledge it has not permeated
quantum chemistry to the same extent.

One important caveat with LOBPCG is its potential numerical
instabilities, which can degrade convergence significantly if care is
not taken. This point was emphasized in \cite{hetmaniuk2006basis},
which describes an appropriately stable procedure to build 
an orthogonal basis of expansion vectors. This relies
crucially on a \texttt{ortho(X,Y)} primitive, whose goal is to
orthogonalize the vectors in $X$ with respect to those in $Y$ and
among themselves. This turns out to be surprisingly hard to perform
reliably in the presence of roundoff errors. The strategy suggested in
\cite{hetmaniuk2006basis} relies on the \texttt{ortho(X,Y)} function
in \cite{stathopoulos2002block}, which is based on a modified singular
value decomposition (SVD). This however can sometimes fail
\cite{duersch2018robust}, and the use of the SVD can become expensive.
For this reason, \cite{duersch2018robust} suggests dropping
ill-conditioned directions.

In this paper, we clarify the origins of numerical instabilities, in
particular due to the ill-conditioned orthogonalizations and reuse of
matrix-vector products. We then present an implementation that uses
exclusively Cholesky decompositions for orthogonalizations. These
decompositions are known to be efficient but unstable, and for this
reason the Cholesky orthogonalization has often been regarded as
impractical; however, it was recently realized that appropriate
stabilizations based on repeated, potentially shifted factorizations,
turn it into a completely reliable orthogonalization algorithm
\cite{fukaya2020shifted}. The resulting LOBPCG algorithm is extremely
stable (converging, if desired, to close to machine precision, with no
degradation in convergence rate), simple (with no dropping of
directions, and with a single parameter controlling the accuracy of
the orthogonalizations), and efficient.

In this paper, we present two opensource implementations of the
{LOBPCG} algorithm. One is \textsc{diaglib}, an opensource,
library written in Fortran 95 that also includes an implementation of
Davidson's method, and that can be obtained at the following address:
\url{https://github.com/Molecolab-Pisa/diaglib}. The library has been
interfaced with the \textsc{CFOUR} suite of quantum chemistry
programs \cite{cfour,matthews2020cfour}, which we have used to test its performance choosing a
variety of test cases. The other is included in the DFTK Julia
software package\cite{DFTKjcon}, supports generalized eigenvalue problems and GPUs,
and is available at \url{https://github.com/JuliaMolSim/DFTK.jl/}.
These implementations are more focused on numerical stability than on
sheer speed, but are overall efficient, thanks to the extensive use of
highly optimized \textsc{BLAS} and \textsc{LAPACK} routines.

We then test the LOBPCG implementation on a selection of quantum chemistry problems, that include full CI, second-order SCF, and CASSCF calculations. For such methods, we compare the performance of LOBPCG to the ones of Davidson's method. 
Our preliminary results show that, while Davidson diagonalization is
more efficient for strongly diagonally dominant problems, as the ones
encountered in Full CI, LOBPCG is a viable alternative, as it exhibits
similar performance and, for difficult cases such in CASSCF, can even
slightly outperform Davidson. 

This paper is organized as follows. In section \ref{sec:theory}, we
present the LOBPCG algorithm. In section \ref{sec:impl}, we discuss in detail the numerical implementation of LOBPCG, analyzing the possible numerical issues and proposing cost-effective, yet robust, solutions. In section \ref{sec:numer}, we test our implementation of a few quantum chemistry applications. Finally, some final remarks are given in section \ref{sec:concl}.

\section{Theory\label{sec:theory}}
In this section, we describe the LOBPCG algorithm in infinite
precision arithmetic, without consideration for numerical stability.
In this paper, we will focus on real symmetric matrices, but
everything of course extends to Hermitian complex matrices.
Let $\mathbf{A} \in \mathbb{R}^{n\times n}$ be a symmetric matrix, for
which we seek $m \ll n$ eigenvalues. Given a set of vectors
$Y = (y_{1}, \dots, y_{p}) \in \mathbb R^{n \times p}$ with $p \ge m$, the Rayleigh-Ritz variational procedure
obtains an approximation ${\rm RR}(Y) \in \mathbb R^{n \times m}$ to the first $m$ eigenpairs by building an
orthonormal basis of $Y$, computing the $p \times p$ matrix representation of $A$ on
that set of vectors, diagonalizing it and taking the first $m$ eigenvectors.

To use this procedure to obtain an iterative algorithm, a way of
constructing $Y$ must be specified. A standard method for doing so is
to use the residuals, defined for a vector $x$ by $r(x) = \mathbf{A} x - \rho(x) x$,
where $\rho(x) = \tfrac{x^{T} \mathbf A x}{x^{T} x}$ is the Rayleigh quotient of
$x$. In ill-conditioned problems, these residuals might not however be
a good search direction, and therefore it is useful to precondition
them according to
\begin{align*}
  w = \mathbf{T}r(x)
\end{align*}
where $T$ is a given preconditioner (for instance, when $\mathbf{A}$ is
diagonally dominant,
$\mathbf{T} = (\mathrm{diag}(\mathbf{A}) - \rho(x) \mathbf{I})^{-1}$).

This choice of search direction results in the block Davidson
algorithm, starting from an initial set of vectors
$X^{[0]} \in \mathbb R^{n \times m}$:
\begin{align*}
  X^{[k+1]} = {\rm RR}(X^{[1]}, \mathbf T^{[1]}R^{[1]} \dots, \mathbf T^{[k]} R^{[k]})
\end{align*}
where $R^{[k]} = \mathbf{A} X^{[k]} - X^{[k]}\Lambda^{[k]}$ is the
residual at the $k$-th iteration, and $\Lambda^{[k]}$ is the diagonal
matrix composed of the Ritz values. Since the expansion subspace is only ever
enlarged, the implementation is standard: each block of vectors added
to the subspace is orthogonalized against the previous vectors and
against itself (although see subtleties of this operation in the next section).

%\textbf{AL: there is some ambiguity here as to whether preconditioning
%is performed on the residuals, the residuals orthogonalized wrt X, or
%something else. What is done in standard Davidson?}

Since the Davidson method performs the Rayleigh-Ritz procedure in the
full convergence history, its computational requirements can increase
quickly. The LOBPCG algorithm instead only keeps the last two iterates:
\begin{align*}
  X^{[k+1]} = {\rm RR}(X^{[k-1]}, X^{[k]}, T^{[k]}(A X^{[k]} - X^{[k]} \Lambda^{[k]})).
\end{align*}
The method is locally optimal (LO) because the Rayleigh-Ritz procedure
optimizes the Rayleigh quotient in the local expansion subspace. It is a
block algorithm (B), and uses a preconditioner (P). Finally, the
intuition of keeping only the previous iterate comes from the
conjugate gradient (CG) algorithm for solving linear systems. This
algorithm can seem like a drastic truncation of the Davidson method.
There is however reason to believe that it can converge asymptotically as quickly as
the full Davidson algorithm, inspired by the optimality in the Krylov
space of the three-terms conjugate gradient algorithm
\cite{knyazev2001toward}.

The convergence properties of this algorithm are sensitive to the gap
between eigenvalues $m$ and $m+1$, which might be small. This is
particularly clear in the case of a simplified version of the LOBPCG
algorithm, the block gradient descent with fixed step (which might be
termed ``BG'', since it is obtained by removing the locally optimal,
preconditioning and conjugate features of LOBPCG), where explicit
convergence rates can be obtained easily \cite{cances2021convergence}.
Accordingly, as is standard, in practice one uses a block size $m$
which is larger than the number of eigenvalues $m_{\rm sought}$
actually sought, and stops the algorithm as soon as the first
$m_{\rm sought}$ eigenvalues are converged. The convergence rate is
then dependent on the gap between eigenvalues $m_{\rm sought}$ and
$m+1$.

\section{Implementation\label{sec:impl}}
\subsection{The LOBPCG algorithm}
When implementing the above algorithm on a computer, we face the
difficulty that the basis
$(X^{[k-1]}, X^{[k]}, T^{[k]}(A X^{[k]} - X^{[k]} \Lambda^{[k]}))$ is
extremely badly conditioned. This is because, as the iteration
progress, $X^{[k]}$ becomes close to $X^{[k-1]}$, and the residual
becomes small. Therefore, if we try to solve the Rayleigh-Ritz problem
as a generalized eigenvalue problem, the results will be inaccurate.
Instead, following \cite{hetmaniuk2006basis}, we construct systematically an orthogonal
basis $(X^{[k]}, W^{[k]}, P^{[k]})$ of the expansion subspace spanned by $(X^{[k-1]}, X^{[k]}, T^{[k]}(A X^{[k]} - X^{[k]} \Lambda^{[k]}))$. The
$P^{[k]}$ is implicitly constructed as the orthogonalization of
$X^{[k]}$ against $X^{[k-1]}$; the $W^{[k]}$ is constructed as the
orthogonalization of $T^{[k]}(A X^{[k]} - X^{[k]} \Lambda^{[k]})$ against
$X^{[k]}$ and $P^{[k]}$.

To obtain these orthogonal bases, we introduce the primitive
\texttt{ortho(X,Y)} which, given a set of orthogonal vectors $Y$,
returns an orthogonal basis of the projection of the vectors in $X$
onto the space orthogonal to $Y$. In infinite precision arithmetic,
this would be given by an orthogonalization of $X - Y Y^{T} X$; in
finite precision arithmetic, care has to be taken, as we will see in
the next section. Given this primitive, the LOBPCG algorithm is given
in Algorithm~\ref{LOBPCG}.
\begin{algorithm}
\caption{\textsc{LOBPCG}\label{LOBPCG}}
\hspace*{\algorithmicindent} \textbf{Input:} Initial guess $X^{[0]}$,
operator $\mathbf{A}$. \\
\hspace*{\algorithmicindent} \textbf{Output:} $X, \lambda$, the lowest
eigenpairs of $\mathbf{A}$.\\
\hspace*{\algorithmicindent} \textbf{Intermediates:} Residuals $R$,
preconditioned residuals $W$, differences $P$, unorthogonalized
quantities (tildes).\\
\begin{algorithmic}[1]
    \State $AX^{[0]} = \mathbf{A} X^{[0]}$
    \State $a^{[0]} = (X^{[0]})^T AX^{[0]}$ \Comment{Initial Rayleigh-Ritz}
    \State solve $a^{[1]}u^{[1]} = u^{[1]} \lambda^{[1]}$
    \State $X^{[1]} = X^{[0]}u^{[1]}$ and $AX^{[1]} = AX^{[0]}u^{[1]}$ 
    \State $R^{[1]} = AX^{[1]} - X^{[1]}\lambda^{[1]}$
    \State $\widetilde{W}^{[1]} = TR^{[1]}$ \Comment{Preconditioned
      residuals...}
    \State $W^{[1]} = $ \textbf{ortho}($\widetilde{W}^{[1]}; X^{[1]}$)
    \Comment{... orthogonalized}
    \State $AW^{[1]} = \mathbf{A}W^{[1]}$
    \State k = 0
    \While{$k \lt k_{\rm max}$}
    \State k = k + 1
    \If {$k = 1$} \Comment{Expansion subspace $V^{[k]}$}
    \State $V^{[k]} = (X^{[k]}, W^{[k]})$
    \Else
    \State $V^{[k]} = (X^{[k]}, W^{[k]}, P^{[k]})$
    \EndIf
    % \State Set $V^{[k]} = (X^{[k]}, W^{[k]}, P^{[k]})$ if $k \ge 2$,
    % $(X^{[k]}, W^{[k]})$ otherwise
    \State $a^{[k]} = (V^{[k]})^T AV^{[k]}$ \Comment{Rayleigh-Ritz in $V^{[k]}$}
    \State solve $a^{[k]}u^{[k]} = \lambda^{[k]} u^{[k]}$, get lowest eigenvectors $u_x^{[k]}$
    \State $X^{[k+1]}$ = $V^{[k]}u_x^{[k]}$ and $AX^{[k+1]}$ = $AV^{[k]}u_x^{[k]}$
    \State $R^{[k]} = AX^{[k]} - X^{[k]}\lambda^{[k]}$
    \State Lock converged eigenvectors, exit if done
    \State Compute $\widetilde u_p^{[k]}$ from unconverged $u_{x}^{[k]}$
    \Comment{Components of $P$ in $V^{[k]}$...}
    \State $u_p^{[k]}$ = \textbf{ortho}($\widetilde{u}_p^{[k]};u_x^{[k]}$)
    \Comment{... orthogonalized}
    \State $P^{[k+1]} = V^{[k]}u_p^{[k]}$ and $AP^{[k+1]} = AV^{[k]}u_p^{[k]}$
    % \State $X^{[k+1]} = X_{\rm new}$ and $AX^{[k+1]} = AX_{\rm new}$
    % \State Set $V^{[k+1]} = (X^{[k+1]},P^{[k+1]})$
    % \State Set $AV^{[k+1]} = (AX^{[k+1]},AP^{[k+1]})$
    \State $\widetilde{W}^{[k+1]} = TR^{[k]}$ \Comment{Preconditioned residuals...}
    \State $W^{[k+1]} = $ \textbf{ortho}($\widetilde{W}^{[k+1]}; (X^{[k+1]}, P^{[k+1]}))$
    \Comment{... orthogonalized}
    \State $AW^{[k+1]} = \mathbf{A}W^{[k+1]}$
% \State $AW^{[k]} = \mathbf{A}W^{[k]}$
% \State Append $W^{[k+1]}$ to $V^{[k+1]}$ and $AW^{[k+1]}$ to $AV^{[k+1]}$
\EndWhile
\end{algorithmic}
\end{algorithm}

%\textbf{TODO revisit after order modification}
\paragraph{Basis selection}
In such algorithm, $X^{[k]}, W^{[k]}, P^{[k]}$ are the $n\times m$ matrices that contain the $m$ desired eigenvectors ($X$) and the corresponding preconditioned residuals ($W$) and increments ($P$), and we denote with a $\sim$ symbol the vectors before orthogonalization. 
The matrices $AX^{[k]}, AW^{[k]}, AP^{[k]}$ contain the results of the application of $\mathbf{A}$ to such vectors. \textsc{LOBPCG} is a \emph{matrix-free} algorithm, i.e., it does not require to assemble and store in memory the matrix $\mathbf{A}$ but just to be able to perform the relevant matrix-vector multiplications. 
% More in detail, an actual matrix-vector multiplication is performed in the beginning with the guess vectors $X^{[0]}$ (line 1) and then once per iteration (line 27). % with the preconditioned residuals $W^{[k]}$, as $AX^{[k]}$ and $AP^{[k]}$ can be obtained indirectly (lines 19,24)
% Two orthogonalizations are performed, one with $n$-sized vectors
% ($\widetilde{W}^{[k]}$ against $V^{[k]}$, line 7, against $X^{[k]}$, and
% line 26), one in the subspace $\widetilde{u_p}^{[k]}$ against $u_x^{[k]}$.
More details on the orthogonalization procedure are given in section \ref{sec:ortho}.
The reduced matrix $a^{[k]} \in \mathbb{R}^{3m\times 3m}$ is diagonalized using standard dense linear algebra routines (in our implementation, we use LAPACK's \textsc{dsyev}).

In the block implementation of \textsc{LOBPCG}, computing $\widetilde{P}^{[k]}$ as $X^{[k]}-X^{[k-1]}$ can become problematic, as these vectors become smaller and smaller when approaching convergence, which can create numerical instabilities. As a more robust alternative, the $\widetilde{P}^{[k]}$ vectors are computed in a different way. Let $u_x^{[k]} \in \mathbb{R}^{3m\times m}$ be the first $m$ eigenvectors of the reduced matrix. 
The eigenvectors have a block structure, that is 
\begin{equation}
    \label{eq:a_block}
    u^{[k]} = \left (\begin{array}{c}
         u_{xx}^{[k]}  \\
         u_{wx}^{[k]}  \\
         u_{px}^{[k]}
    \end{array}
    \right )
\end{equation}
where each block is a $m\times m$ square matrix. 
The new eigenvectors are computed as (line 19)
\begin{equation}
    \label{eq:neweigv}
    X^{[k+1]} = X^{[k]}u_{xx}^{[k]} + W^{[k]}u_{wx}^{[k]} + P^{[k]}u_{px}^{[k]}
\end{equation}
To compute the $P^{[k+1]}$ vectors, we first get the expansion
coefficients $\widetilde{u_p}^{[k+1]}$ of $X^{[k+1]}-X^{[k]}$ in
$V^{[k]}$, which are obtained by subtracting the identity matrix from
the unconverged components of $u_{xx}^{[k]}$. Then, we orthogonalize them against $u_x$,
and use them to compute the new $P^{[k+1]}$ vectors. Note
that we assemble the $P$ vectors only corresponding to the active eigenvectors,
i.e., the ones that have not yet converged; for this reason, it is
important to perform such operation before orthogonalizing
$\widetilde{u_{p}}^{[k+1]}$.

\paragraph{Reuse of applications}
Apart from the choice of a basis and its orthogonalization, a
numerically sensitive point is the reuse of the applications of
$\mathbf A$. Since this is a potentially costly operation, it is not
feasible to recompute for instance $\mathbf A P^{[k]}$ before the
Rayleigh-Ritz procedure; instead, we use the fact that $P^{[k]}$ is
built as a linear combination of other vectors, on which we know the
application of $\mathbf A$. If this is done naively however this can
result in a large error. This is because, in general, if $\mathbf A V$
is known to some precision $\varepsilon$, then $(\mathbf A V) u$ will
be an approximation of $\mathbf A (Vu)$ with a precision of the order
of $\|u\| \varepsilon$.

Consider the problem of computing $\mathbf {A} P^{[k]}$, line 24 of the
algorithm. In exact arithmetic, we could compute
$\widetilde{P}^{[k]} = V^{[k]} \widetilde u_{p}^{[k]}$, compute
$\mathbf A \widetilde{P}^{[k+1]} = (\mathbf{A} V^{[k]}) \widetilde u_{p}^{[k]}$,
orthogonalize $\widetilde{P}^{[k+1]}$ against $X^{[k+1]}$ and update
$\mathbf{A} \widetilde P^{[k+1]}$ accordingly, etc. This however
amounts to obtaining $\mathbf{A} P^{[k+1]}$ by right-multiplying
$\mathbf{A}V^{[k]}$ with a sequence of potentially ill-conditioned
(and therefore of large norm) matrices, which incurs a large error on
$\mathbf{A} P^{[k+1]}$. Instead, we obtain directly the expansion
coefficients $u_{p}^{[k]}$ of $P^{[k+1]}$ on $V^{[k]}$, and
obtain $\mathbf{A} P^{[k+1]}$ as $(\mathbf{A} V^{[k]}) u_{p}^{[k]}$. Since
both $P^{[k+1]}$ and $V^{[k]}$ are orthogonal, so is $u_{p}^{[k]}$, and
therefore no precision is lost in the update
$AP^{[k+1]} = AV^{[k]}u_p^{[k]}$. The same is true for the update
$AX^{[k+1]} = AV^{[k]}u_x^{[k]}$ (line 19).

% In a naive implementation, $u_p$ would contain the
% coefficients of $P^{[k]}$ in the expansion subspace, and could be 
% very ill conditioned, %Since here $u_p$ results from the orthogonalization
% %of potentially very badly conditioned vectors, 
% and therefore its norm can be large,
% resulting in a large error in $\mathbf A P^{[k]}$. This in turn can
% induce large errors in the Rayleigh-Ritz procedure. Rather, we make
% sure to only use the rule $\mathbf A (V u_p)= (\mathbf A V) u_p$ when $u_p$
% is an orthogonal matrix. In the specific case, we first orthogonalize $u_p$ to $u_x$,
% and then proceed safely to compute $P = Vu_p$ and $AP = AVu_p$.
% The safe procedure described in this paragraph is further used to compute both
% $AX^{[k+1]}$ (line 19) and $AP^{[k+1]}$ (line 27).

\paragraph{Locking}
Another crucial aspect of an efficient and stable implementation concerns the treatment of converged eigenvectors. In our implementation, we freeze the first $m_{\rm conv}$ consecutive eigenvectors, which means that we only compute $m_{\rm act} = m - m_{\rm conv}$ new residuals, $W$ and $P$ vectors. The converged eigenvectors are kept into $X$, to enforce the orthogonality of the active search subspace. This means that the reduced matrix and the $V$ subspace dimensions are $m + 2m_{\rm act}$, and that only $m_{\rm act}$ matrix-vector multiplications are performed at each iteration, combining thus stability and efficiency.

\subsection{A robust and stable \texttt{ortho(X,Y)}
  procedure\label{sec:ortho}}
One of the most crucial steps in \textsc{LOBPCG} is the
orthogonalization of a set of vectors against a given set, and its
subsequent orthonormalization. We first tackle the \texttt{ortho(X)}
routine, which orthogonalizes a set of vectors.

\paragraph{The \texttt{ortho(X)} procedure}
The gold standard for orthogonalizing a set of vectors is to compute the (thin) singular value
decomposition of $X$ and then take the left singular vectors. A
slightly less expensive, yet very stable alternative, is to use the QR decomposition of X, which in our tests performs equivalently well.
Another good option is the modified Gram-Schmidt algorithm. However,
these algorithms can all become expensive, especially if a large
number of eigenvalues are sought.

An alternative and cheaper strategy is to compute the Cholesky
decomposition of the overlap matrix %\textbf{AL: one advantage of Cholesky
  %(and the main motivation for me to look at it at the time) is that
  %it parallelizes much better, esp. in the multi-node (MPI) context.
  %Do we want to say anything about that here, or do people in
  %chemistry still freak out about multi-node?}
\begin{equation}
    \label{eq:CDOrtho}
    X^{T} X = M = LL^T
  \end{equation}
%  \textbf{AL: I'm not sure about calling $M$ the metric, as it feels a
%  bit generic. How about overlap matrix?}
The orthogonal vectors can then be obtained by solving the triangular linear system
\begin{equation}
    \label{eq:CDOrtho2}
    {\texttt{ortho}}(X) L^T = X.
  \end{equation}
This is often more efficient, as it allows for greater
parallelization and full use of BLAS3 routines.
  
  This procedure works in infinite precision, but has two issues in
  finite precision. First, even after a first Cholesky orthogonalization, the vectors
  can fail to be orthogonal. Second, the Cholesky decomposition can fail,
  because $M$ may not be positive-definite to machine precision;
  this happens when the conditioning of $X$ is larger than the square
  root of the inverse machine epsilon, about $10^{8}$ in double
  precision arithmetic.

  Fortunately, there is a simple fix to the first problem:
  orthogonalize twice. This has been established to produce vectors
  orthogonal to machine precision \cite{yamamoto2015roundoff}. For the
  second problem, following \cite{fukaya2020shifted}, we level shift
  the metric before its Cholesky decomposition by adding a small
  constant to its diagonal. Such a constant can be chosen very small
  (in our implementation, we start from 100 times the norm of $X$ times
  the machine precision) and, if the decomposition still fails,
  increased until the Cholesky decomposition is successful. In our
  tests, the shifted decomposition never failed, and therefore, with at most 4
  Cholesky orthogonalizations (one failed unshifted, one shifted, then
  two unshifted) we are guaranteed to obtain vectors that are
  orthogonal to machine precision. In practice, often much less than
  this is needed -- failures of the first Cholesky orthogonalization
  have been observed only exceptionally.
  A pseudo-code for the \texttt{ortho} procedure is given in algorithm \ref{alg:ortho_cd}.

\begin{algorithm}
\caption{\texttt{ortho}$(\widetilde X)$: Orthonormalize a set of vectors $\widetilde{X}$ using the Cholesky decomposition of the overlap with iterative refinement.}
\hspace*{\algorithmicindent} \textbf{Input:}  non orthogonal vectors $\widetilde{X}$, threshold $\tau_{\rm ortho}$. \\
\hspace*{\algorithmicindent} \textbf{Output:} $X$, orthonormal vectors.%\\
% \hspace*{\algorithmicindent} \textbf{Intermediates:} Overlap matrix $M$.
% \\
\label{alg:ortho_cd}
%  \begin{algorithmic}[1]
%  \State k = 0
%  \While{$k \lt k_{\rm max}$}
%      \State $k = k + 1$
%      \State $M = \widetilde{X}^T\widetilde{X}$
%      \State Attempt Cholesky factorization $M = LL^T$
%      \If{fail}
%          \State Add $\alpha \epsilon$ to the diagonal of $M$ until successful
%          % \State $\alpha = 100$
%          % \State l = 0
%          % \While{$l \lt l_{\rm max}$}
%          %     \State $l = l + 1$
%          %     \State $\widetilde{M} = M + \alpha \epsilon \|U\| Id$
%          %     \State $\widetilde M = LL^T$
%          %     \If{ok}
%          %         \State exit
%          %     \EndIf
%          % \EndWhile
%      \EndIf
%      \State $X = \widetilde{X} L^{-T}$
%      \State $M = X^TX$
%      \If{$\|M - Id\| < \tau_{\rm ort}$}
%          exit
%      \Else
%         \State $\widetilde{X} = X$
%      \EndIf
%  \EndWhile
% \end{algorithmic}

 \begin{algorithmic}[1]
   \State $X = \widetilde X$
 \While{$\|X^{T} X - Id\| > \tau_{\rm ort}$}
     \State $M = {X}^T{X}$
     \State Attempt Cholesky factorization $M = LL^T$
     \If{fail}
         \State Add $\alpha \epsilon \|X\|$ to the diagonal of $M$ until successful
         % \State $\alpha = 100$
         % \State l = 0
         % \While{$l \lt l_{\rm max}$}
         %     \State $l = l + 1$
         %     \State $\widetilde{M} = M + \alpha \epsilon \|U\| Id$
         %     \State $\widetilde M = LL^T$
         %     \If{ok}
         %         \State exit
         %     \EndIf
         % \EndWhile
     \EndIf
     \State $X = {X} L^{-T}$
 \EndWhile
 \end{algorithmic}

 \end{algorithm}
%  \textbf{TODO give some stats here? And do we describe
%    the norm estimation trick to avoid repeating the
%    orthogonalizations?}

 The algorithm as given is somewhat wasteful, as in the common case
 where only one or two successful Cholesky factorizations are needed
 it recomputes the overlap to check for termination. This can
 potentially be alleviated by computing a cheap estimation to the norm
 of $L^{-T}$: if this is moderate, then the new vectors are orthogonal
 to a good accuracy, and a new round is unnecessary.

  \paragraph{The \texttt{ortho(X,Y)} procedure}
  Using the previous orthogonalization algorithm, we could implement the
  ${\tt ortho}(X,Y)$ as ${\tt ortho}(X - Y Y^{T} X)$. This is however
  numerically unstable: if $X - Y Y^{T} X$ is of order $\delta$
  (because $X$ was almost in the range of $Y$), then the
  orthogonalization above will multiply it by a factor of order $1/\delta$,
  meaning that $Y^{T} {\rm ortho}(X - Y Y^{T} X)$ will be of order
  $\varepsilon/\delta$, where $\varepsilon$ is the machine epsilon,
  and the vectors will not be sufficiently orthogonal to $Y$. To avoid
  this, we use a loop: first project out $Y$, then orthogonalize, iteratively until convergence. In practice, two steps
  are usually enough to achieve
  convergence. % we needed \textbf{AL: give some ideas of how many steps
  %are needed in practice}
  The algorithm for the \texttt{ortho(X,Y)} procedure is given in algorithm \ref{alg:ortho}.
 \begin{algorithm}
 \caption{\texttt{ortho}$(\widetilde X,Y)$: given a set of orthonormal vectors $Y$ and a set of vectors $\widetilde{X}$, orthogonalize $\widetilde{X}$ to $Y$ and orthonormalize $\widetilde{X}$. }
 \hspace*{\algorithmicindent} \textbf{Input:}  orthonormal vectors $Y$, non orthogonal vectors $\widetilde{X}$, threshold $\tau_{\rm ortho}$. \\
\hspace*{\algorithmicindent} \textbf{Output:} $X$, orthonormal vectors also orthogonal to $Y$.
 \label{alg:ortho}
 \begin{algorithmic}[1]
 % \State $X$ = \textbf{ortho}($\widetilde{X}$)
 \State $X$ = $\widetilde{X}$
 \While{$\|Y^{T}X\| > \tau_{\rm ortho}$}
     \State $X = X - YY^TX$
     \State $X =$ \textbf{ortho}(${X}$)
 \EndWhile
 \end{algorithmic}
 \end{algorithm}
%  \begin{algorithm}
%  \caption{\texttt{ortho}$(\widetilde X,Y)$: given a set of orthonormal vectors $Y$ and a set of vectors $\widetilde{X}$, orthogonalize $\widetilde{X}$ to $Y$ and orthonormalize $\widetilde{X}$. }
%  \hspace*{\algorithmicindent} \textbf{Input:}  orthonormal vectors $X$, non orthogonal vectors $\widetilde{X}$, threshold $\tau_{\rm ortho}$. \\
% \hspace*{\algorithmicindent} \textbf{Output:} $X$, orthonormal vectors also orthogonal to $Y$.\\
% \hspace*{\algorithmicindent} \textbf{Intermediates:} Overlap matrix $M$.
%  \label{alg:ortho}
%  \begin{algorithmic}[1]
%  \State k = 0
%  \State X = \textbf{ortho}($\widetilde{X}$)
%  \While{$k \lt k_{\rm max}$}
%      \State $k = k + 1$
%      \State $\widetilde{X} = X - YY^TX$
%      \State X = \textbf{ortho}($\widetilde{X}$)
%      \State $M = X^T X$
%      \If{$\|M\| < \tau_{\rm ort}$}
%          \State exit
%      \EndIf
%  \EndWhile
%  \end{algorithmic}
%  \end{algorithm}

 Similar to before, this algorithm is relatively wasteful in the common case
 where one or two passes are enough, because it recomputes $Y^{T}X$ to
 check for termination. This can be remedied by monitoring the growth
 factor of \texttt{ortho}$(X)$ (the maximum amplification of errors in
 $X$ caused by the \texttt{ortho} routine, and therefore a measure of
 the lack of $Y$-orthogonality after one iteration), and exiting the
 loop when that number is moderate.

\subsection{Generalized eigenvalue problems}
Missing from this algorithm is a discussion of generalized eigenvalue
problems, simply because they are not often encountered in quantum
chemistry. In the generalized eigenvalue problem, one solves
$\mathbf{A} x = \lambda \mathbf{B} x$, where $\mathbf{B}$ is a
symmetric positive definite matrix. Eigenvectors are orthogonal with
respect to the modified inner product
$\langle x, y \rangle_{\mathbf B} = x^{T} \mathbf{B} x$. The
theoretical LOBPCG algorithm is unchanged except for the fact that
residuals are now ${\mathbf A} x - \lambda {\mathbf B} x$, and that
all orthogonalizations are with respect to the modified inner product.

Our practical algorithm has to be modified by keeping a
$\mathbf{B}$-orthogonal basis $V$. This could be done by using the
$\textbf{B}$ inner product in the orthogonalization of $W$, line 26 of
the algorithm above, and maintaining the values of $\mathbf{B} V$
along the iterations. However, a naive implementation of this step
requires either multiple applications of $\mathbf{B}$ (twice per
iteration, on $\widetilde W^{[k+1]}$ and on $W^{[k+1]}$), or
potentially unsafe reuses of applications of $\mathbf{B}$. As a
compromise, a good option is to use the intermediate quantity
\begin{align*}
  \widehat W^{[k+1]} = \mathtt{ortho}(\widetilde W^{[k+1]}, (\mathbf{B} X^{[k+1]}, \mathbf{B} P^{[k+1]})),
\end{align*}
which is $\mathbf B$-orthogonal to $X^{[k+1]}$ and $P^{[k+1]}$, but
whose vectors are only orthogonal (and not $\mathbf B$-orthogonal) to
each other. This set of vectors $\widehat W^{[k+1]}$ is however well
conditioned (with respect to the $\mathbf B$-inner product). We can
therefore $\mathbf{B}$-orthogonalize it to compute
$W^{[k+1]} = \widehat W^{[k+1]} L^{-T}$ with $L$ a well-conditioned
matrix; it is then safe to re-use the $\mathbf B$ application as
$\mathbf{B} W^{[k+1]} = (\mathbf{B} \widehat W^{[k+1]}) L^{-T}$. This
appeared to perform very well in our tests, even if $\mathbf{B}$ itself was not
well-conditioned. If more stability is needed, reuses of $\mathbf{B}$
applications appear necessary.

\section{Numerical experiments\label{sec:numer}}
To test our implementation of LOBPCG, and to compare its performance
with respect to the block-Davidson method, we interfaced the
\textsc{diaglib} library with the \textsc{CFOUR} quantum chemistry
package \cite{cfour,matthews2020cfour}. As typical problems where an iterative procedure to compute one or a few eigenvectors is required, we selected three different test cases, coming from full CI Hamiltonian, quadratically convergent self-consistent field (SCF), and quadratically convergent complete active space self-consistent field (CASSCF).
In all the calculations, we use a threshold of $10^{-14}$ for the \texttt{ortho(X,Y)} and \texttt{ortho(X)} procedures. Together with the convergence threshold for the eigenvectors, this is the only parameters that control the LOBPCG calculation. For Davidson, we use a subspace dimension of 25, i.e., we keep in memory up to 25 vectors per eigenvector in the history. For both the algorithms we exploit a locking procedure for the converged eigenvectors. 

\subsection{Full CI calculations}
We compute the first few totalsymmetric %\textbf{AL: totalsymmetric $\rightarrow$ totally symmetric?} %we use totalsymmetric in computational chemistry
electronic states of water at the Full CI level of theory, using a determinant CI direct implementation. The Full CI Hamiltonian is sparse and diagonally dominant, but extremely large, and thus provides a good testcase for well-behaved, large, sparse systems. Furthermore, the iterative solution of (Full) CI problems is a quite common task in quantum chemistry, as it is encountered in CASCI/CASSCF, and truncated CI (including for excited states at the CI singles level of theory). We use Pople's 6-31G$^*$ basis set,\cite{hehre1972self} and perform both all-electron and frozen-core calculations, correlating thus 10 electrons in 18 orbitals (18~360~640 determinants) or 8 electrons in 10 orbitals (1~416~732 determinants). We seek 10, 20, or 50 eigenpairs. Convergence is achieved when the root-mean-square norm of the residual is smaller than $10^{-9}$, and its maximum absolute value is smaller than $10^{-8}$. For LOBPCG, we seek 5 additional eigenpairs, as numerical tests proved that this improves convergence and, despite the additional matrix-vector products required, improves overall performance. Note that we do not check that the additional eigenvalues are converged, as they are only used to increase the expansion subspace. No additional eigenpairs are sought for Davidson, as this choice showed the overall best performance. A brief description of the process that led us to these choices is reported in the Supporting Information.

\begin{figure}
    \centering
    \includegraphics[width=0.45\textwidth]{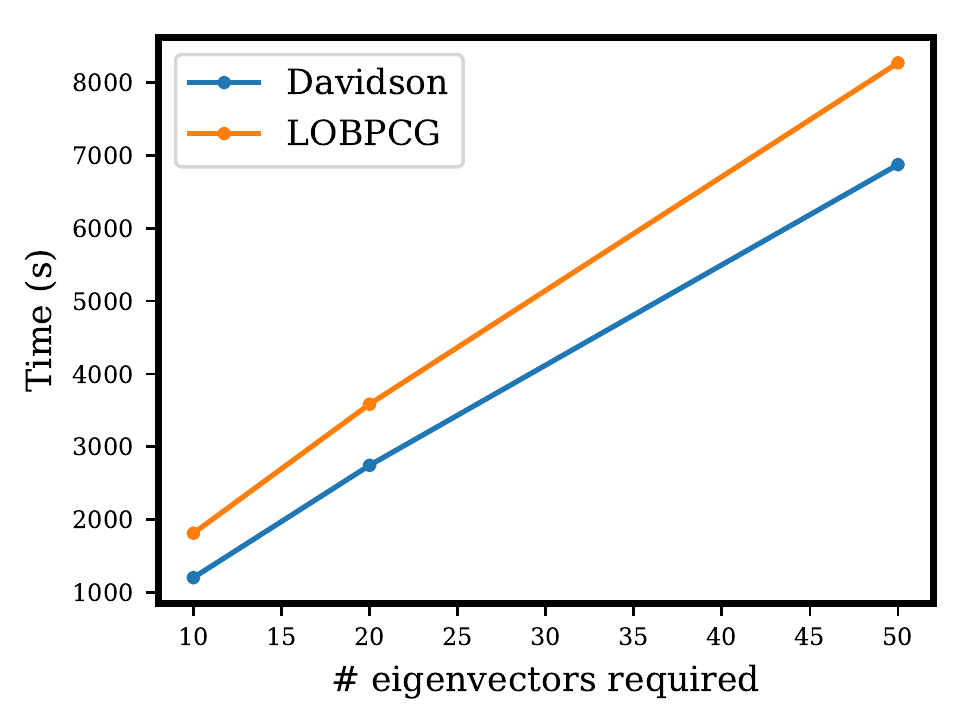}
    \includegraphics[width=0.45\textwidth]{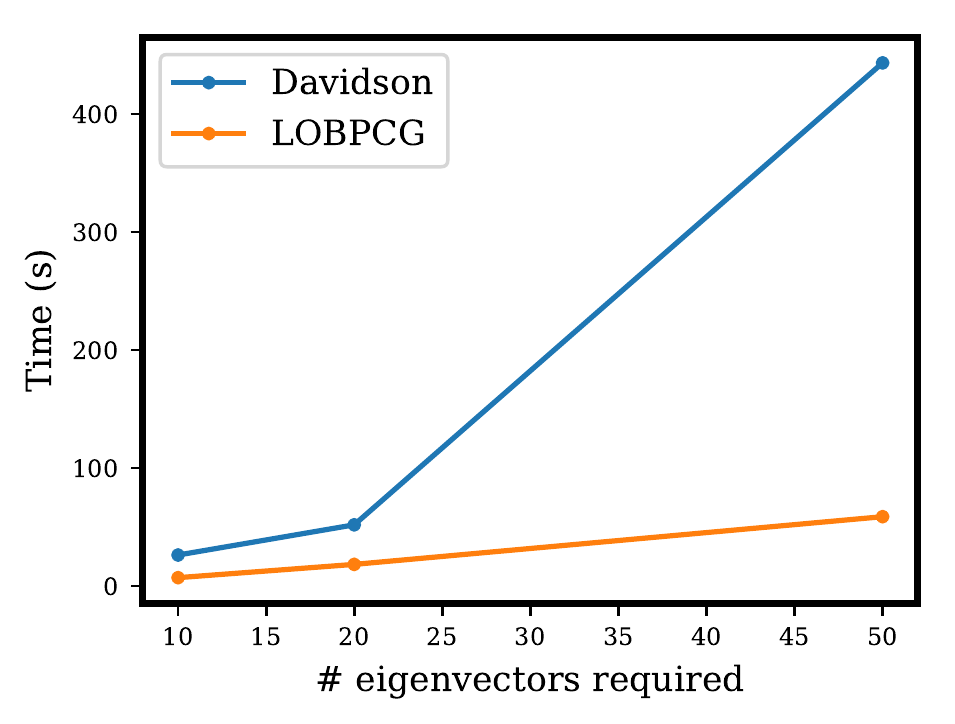}
    \caption{All-electron Full CI calculations for water using Davidson or LOBPCG. Total timings (left panel) and cumulative time for the Rayleigh-Ritz and orthogonalization procedures (right panel).}
    \label{fig:fci_ae}
\end{figure}
\begin{figure}
    \centering
    \includegraphics[width=0.45\textwidth]{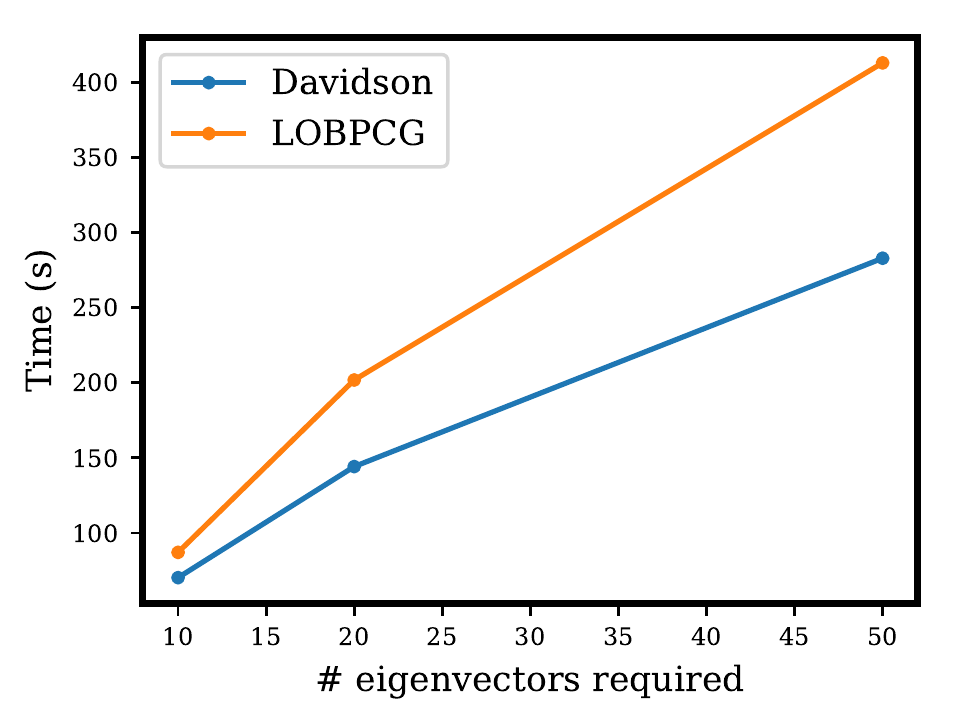}
    \includegraphics[width=0.45\textwidth]{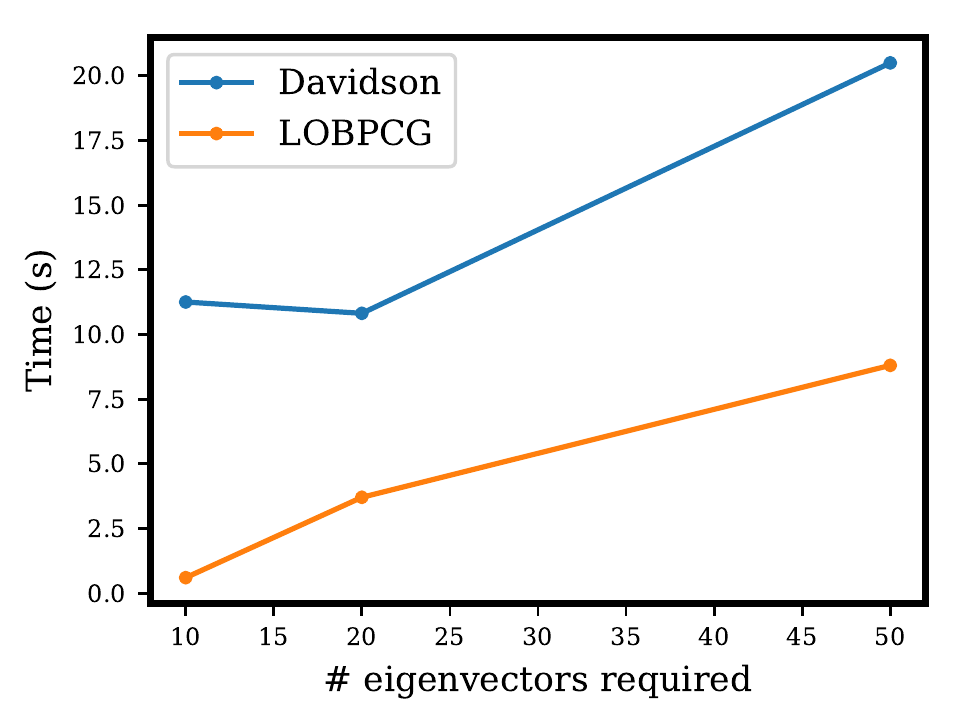}
    \caption{Frozen-core Full CI calculations for water using Davidson or LOBPCG. Total timings (left panel) and cumulative time for the Rayleigh-Ritz and orthogonalization procedures (right panel).}
    \label{fig:fci_fc}
\end{figure}
It comes as no surprise that Davidson outperforms LOBPCG for Full-CI calculations. As the full-CI Hamiltonian is strongly diagonally dominant for closed-shell systems, it is an ideal scenario for Davidson, an algorithm originally conceived for exactly this problem. Nevertheless, the performance of LOBPCG are comparable, the latter algorithm being about 20-30\% slower than the former. It is interesting to note that, while for LOBPCG the cost of the Rayleigh-Ritz and orthogonalization procedures is overall negligible with respect to the cost of computing matrix-vector multiplications, this is not the case for Davidson. Keeping up to 25 vectors in the history comes with a cost, that can be clearly seen in the right panels of figures \ref{fig:fci_ae}, \ref{fig:fci_fc}. On the other hand, the larger subspace used in Davidson's method allows for faster convergence, which is achieved in 28, 37, and 24 iterations for the calculation seeking 10, 20, and 50 eigenpairs, respectively, both for the all-electron and the forzen-core calculations. This has to be compared with 26, 41, and 45 iterations for LOBPCG, again, for both sets of calculations. On the other hand, the long history is also a limitation for Davidson, as the amount of memory required to perform a calculation can become very high. As an example, the largest calculation performed (all-electron, 50 states) required about 356 GB of memory for Davidson, to be compared with 55 for LOBPCG. Using a smaller subspace dimension in Davidson is of course possible, but such a size must be chosen with some care. To better illustrate this point, we repeated the Davidson calculations using a maximum of 10 points in the history: no calculation fully converged within 100 iterations, with 1, 3, and 2 non converged roots. Therefore, while Davidson is optimal for Full-CI calculations if enough memory is available to use a large expansion subspace, LOBPCG can be seen as a competitive alternative when this is not the case.

\subsection{Quadratically convergent SCF calculations}
In quadratically convergent implementations of the self-consistent
field, the Hartree-Fock wavefunction is optimized by using a
second-order method. In \text{CFOUR}, this is done using an efficient
numerical realization of the Levenberg-Marquardt
method \cite{fletcher2013practical}, known as norm-extended
optimization, where the step is computed from the lowest eigenpair of
the (augmented) orbital-rotation Hessian \cite{nottoli2021black}. The same matrix is used in response calculations and for the analysis of the stability of the Hartree-Fock wavefunction, which requires again to compute a few eigevalues and eigenvectors of the orbital-rotation Hessian. For closed-shell systems, such a matrix is dense, but typically diagonally dominant. While computing and storing in memory the full Hessian is possible, such a task is expensive and exhibits a steep scaling in computational cost with respect to the system's size, as assembling it requires a costly partial integral transformation. Direct implementations are therefore usually preferred. 
%To compare Davidson to LOBPCG, we compute the first few eigenpairs (up to 10) of the Hessian after convergence of the SCF procedure for a medium-sized molecule, cumarin, using Pople's 6-31G$^*$ basis set. To provide an example on a somewhat more challenging case, we further repeat the calculations, but at the beginning of the SCF calculations, that is, using orbitals computed with an extended H\"uckel guess. Far from convergence, the Hessian has many negative eigenvalues and is not guaranteed to be as diagonally dominant as with fully converged orbitals. In the following calculations, convergence is achieved when the root-mean-square norm of the residual is smaller than $10^{-7}$, and its maximum absolute value is smaller than $10^{-6}$. Such threshold are adequate for stability analysis, but also for the first steps of a second-order optimization procedure. For both Davidson and LOBPCG, to improve convergence, we seek twice the eigenpairs required, and stop the calculation when just the required ones are converged. The size of the subspace for Davidson is 25.
To compare Davidson to LOBPCG, we compute the first few eigenpairs (up to 10) of the Hessian after convergence of the SCF procedure for a transition metal complex, FeC(CO)$_3$, using Dunning's \textit{cc}-pVDZ basis set. To provide an example on a somewhat more challenging case, we further repeat the calculations, but at the beginning of the SCF calculations, that is, using orbitals computed with an extended H\"uckel guess. Far from convergence, the Hessian has many negative eigenvalues and is not guaranteed to be as diagonally dominant as with fully converged orbitals. In the following calculations, convergence is achieved when the root-mean-square norm of the residual is smaller than $10^{-7}$, and its maximum absolute value is smaller than $10^{-6}$. Such threshold are adequate for stability analysis, but also for the first steps of a second-order optimization procedure. For both Davidson and LOBPCG, to improve convergence, we seek twice the eigenpairs required, and stop the calculation when just the required ones are converged. The size of the subspace for Davidson is 25.
\begin{figure}
    \centering
    \includegraphics[width=0.45\textwidth]{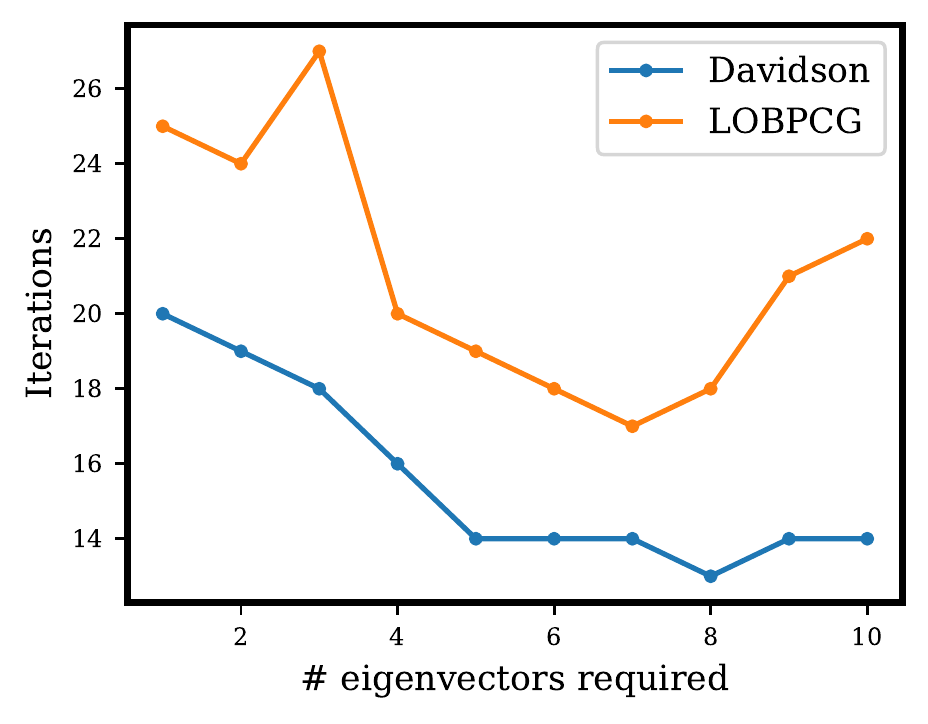}
    \includegraphics[width=0.45\textwidth]{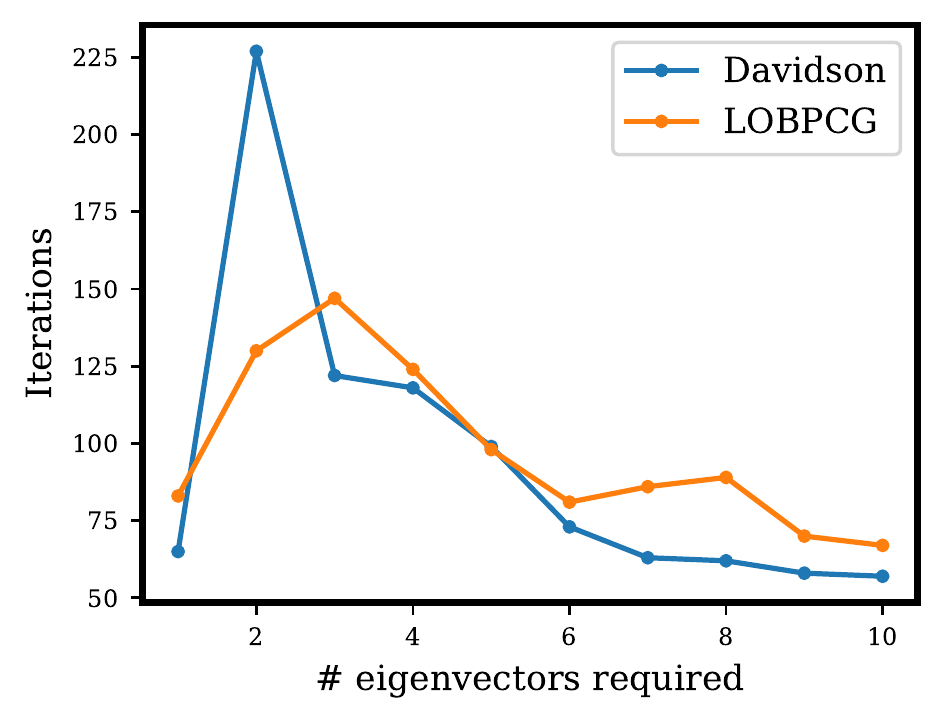}
    \caption{Orbital rotation Hessian diagonalization at SCF convergence (left panel) and with extended-H\"uckel guess orbitals (right panel). The number of iterations required to achieve convergence is reported as a function of the number of seeked eigenpairs.}
    \label{fig:qcscf}
\end{figure}
For the calculations at SCF convergence, where the Hessian is strongly diagonally dominant (figure \ref{fig:qcscf}, left panel) Davidson outperforms again LOBPCG. On the other hand, the less well-behaved case (figure \ref{fig:qcscf}, right panel) shows a different picture. LOBPCG and Davidson exhibit a very similar behavior, with LOBPCG even outperforming Davidson in a few cases. As orbital rotation Hessians are hardly very large matrices (for the case reported here, the size is 3849), using large expansion subspaces in Davidson can probably further improve convergence, but for difficult cases, LOBPCG can be a valid alternative.

\begin{figure}
    \centering
    \includegraphics[width=0.45\textwidth]{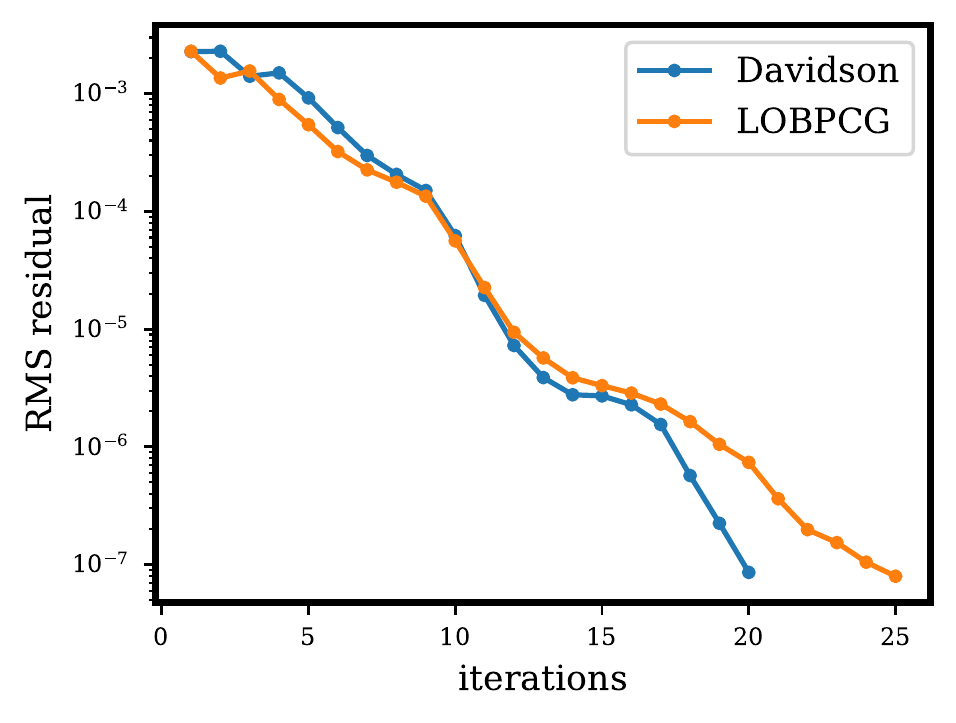}
    \includegraphics[width=0.45\textwidth]{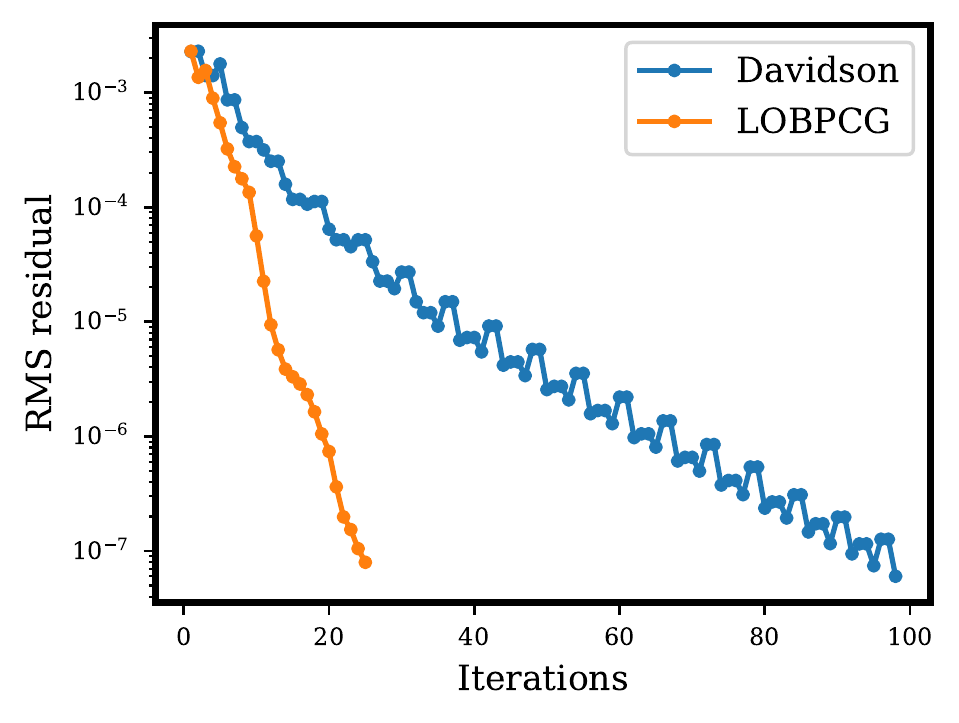}
    \caption{Root mean square (RMS) of the residual along the iterations for the diagonalization of the orbital rotation Hessian at SCF convergence seeking for one eigenpair. The subspace for the Davidson algorithm is of size 25 (left panel) and 3 (right panel). }
    \label{fig:rms_it}
\end{figure}

We report in figure \ref{fig:rms_it} (left panel), the root mean square residual as a function of the iterations for the computation of the first eigenpair of the orbital rotation Hessian. As expected the behaviors between the two algorithms are analogous.

We would like here to underline the fact that LOBPCG manages to behave similarly to Davidson even despite the very small size of the subspace used for the Rayeligh-Ritz procedure. To show how remarkable this is, we report in figure \ref{fig:rms_it} (right panel) a comparison between LOBPCG and Davidson, where for the latter we use a three dimensional subspace -- that is, the same dimension used in LOBPCG. Davidson eventually manages to converge, however, it requires a much larger number of iterations. While this is purely an academic example, as such small expansion subspaces are never used in practice, it testifies to the effectiveness of the LOBPCG 3-terms sequence. 
%Instead, if were to enforce the Davidson subspace to be equal to 3 -- thus requiring a similar amount of memory as for LOBPCG -- we note that the number of iterations rises steeply. Therefore, within this scenario LOBPCG outperforms Davidson as reported in figure \ref{fig:rms_it_dav3}.

\subsection{CASSCF calculations}
\begin{figure}
    \centering
    \includegraphics[width=0.45\textwidth]{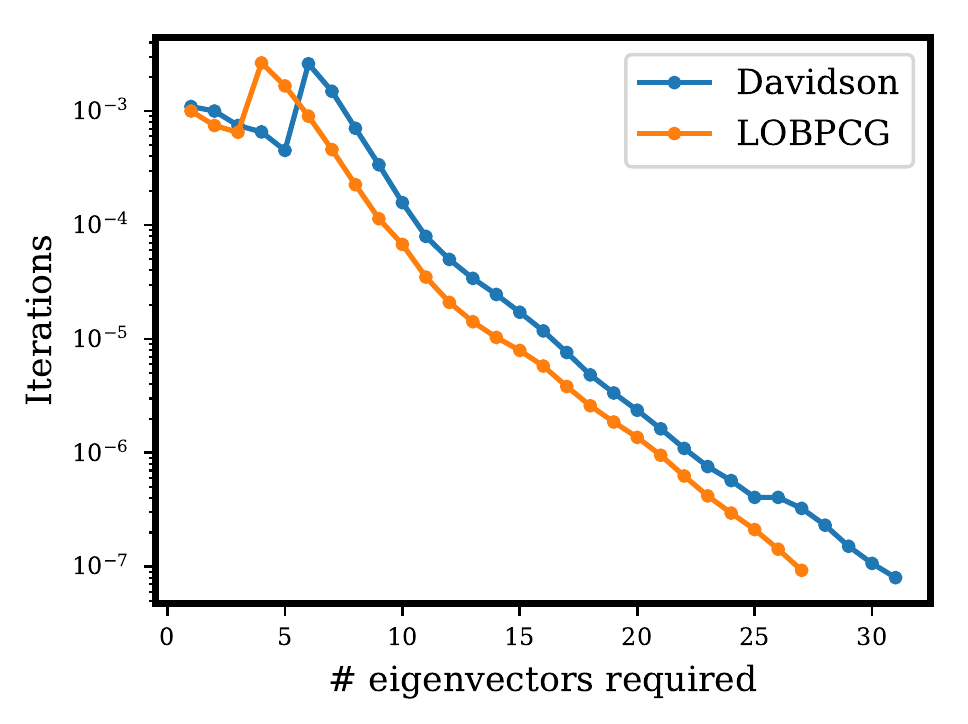}
    \includegraphics[width=0.45\textwidth]{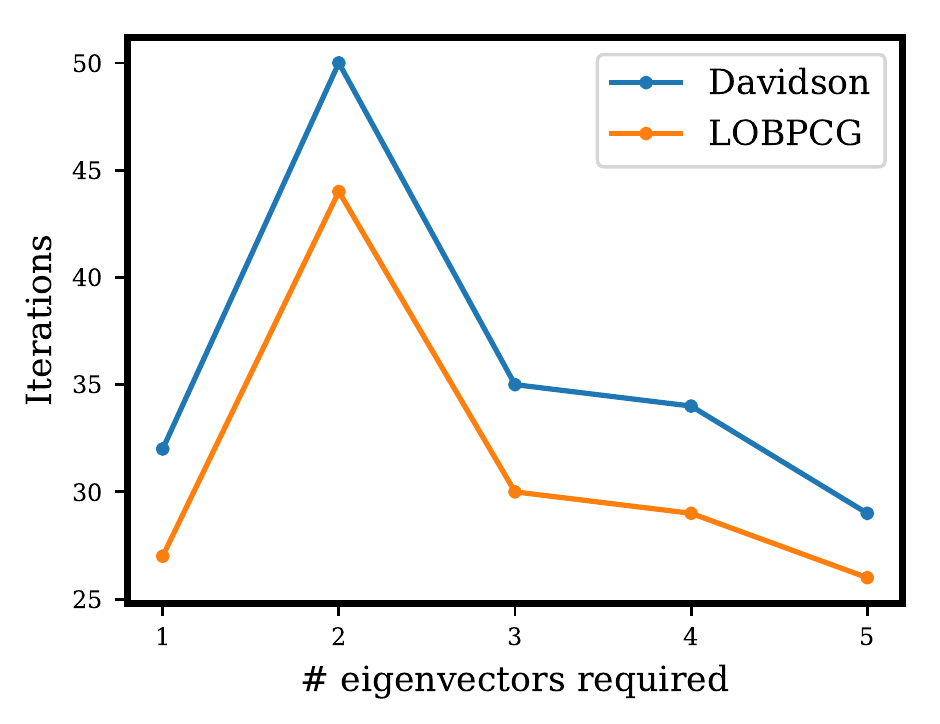}
    \caption{Root mean square (RMS) of the residual along the iterations seeking for one eigenpair (left panel) and number of iterations to achieve convergence as a function of the number of seeked eigenpairs (right panel) for the diagonalization of the CASSCF augmented Hessian at the first macroiteration.}
    \label{fig:casconv}
\end{figure}

CASSCF calculations can be very challenging from a numerical point of
view, which makes second-order methods particularly attractive
\cite{werner1985second,werner1980quadratically,jensen1984direct,jensen1986direct}.
In \textsc{CFOUR}, the same technique used for Hartree-Fock, namely,
the norm-extended optimization algorithm is used. The (augmented)
Hessian in CASSCF is made by a dense, typically quite ill-conditioned,
medium-sized block for the orbital optimization, and a large, sparse,
usually diagonally dominant block that corresponds to the Hamiltonian
in the CAS space. Even for well-behaved systems, computing the NEO
step, which in turn requires computing the first eigenpair of the
augmented Hessian, can be challenging. To illustrate this, we report
calculations on niacin (vitamin B3), a small conjugated organic
molecule. We correlated all the $\pi$ electrons, resulting in a
CAS(6,6) calculation, and we employ Pople's 6-31G* basis
set \cite{hehre1972self}. Symmetry broken unrestricted natural
orbitals are used as a guess \cite{pulay1988uhf,toth2020comparison}. The system is very well behaved, and convergence of the wavefunction is achieved in just 3 second-order iterations. 
Nevertheless, the iterative calculation of the step, i.e., the augmented hessian lowest eigenvector, can still be challenging. 
We report in the left figure \ref{fig:casconv} the convergence pattern for Davidson and LOBPCG at the first second-order iteration.
LOBPCG outperforms Davidson at every second-order iterations, 27, 25, and 23 iterations at the first, second, and third second-order step, to be compared with 32, 28, and 31 iterations for Davidson.
Also when more than one eigenpair is required, as it is the case for state-specific excited state calculations, LOBPCG keeps being the best performing algorithm, as reported in the right panel of figure \ref{fig:casconv}.
This behavior is not surprising, as the CASSCF augmented Hessian is
not diagonally dominant, and is consistent with what observed for
second-order SCF calculations starting with extended H\"uckel orbitals
Davidson's performance may be improved by increasing the size of the expansion subspace. This is, however, not the best option, as for large active spaces the size of the Hessian can become very large -- comparable with the sizes reported for Full CI calculations. Using very large expansion subspaces is therefore very demanding in terms of memory, and can make the orthogonalization expensive. LOBPCG seems therefore a more robust choice for this specific problem.

\subsection{Preconditioning}
To improve the convergence of both LOBPCG and Davidson one can devise different types of preconditioners. In most quantum chemistry applications, computing and storing in memory the matrix of which one seeks one or a few eigenpairs is prohibitively expensive, which forces the choice of a Jacobi (diagonal) preconditioner in most cases. This is definitely the case for the full CI and CASSCF examples showed in the previous section. However, for both second-order SCF and CASSCF, we have implemented the explicit construction of the orbital rotation Hessian, mainly as a debug option, which allows us to perform a few numerical experiments. 
% ($
%\mathbf{T}$).
%\begin{align}\label{eq:precond}
%  \mathbf{T}w = r(x)
%\end{align}
In particular, we compare three possible choices, and focus on the CASSCF orbital-rotation Hessian as a test case, as it is notoriously ill-conditioned and therefore we expect that a more advanced preconditioning strategy may be particularly beneficial. For these examples, we only seek to compute one eigenpair. The first preconditioner that we test is, as in the previous sections, appoximates the matrix to its diagonal. The second one, which will be here addressed as tridiagonal, improves upon the diagonal approximation by also including the upper and lower diagonal elements. These two options should perform particularly well in diagonally-dominant matrices. As a third choice, we propose a sparse approximation $M$ to the matrix $A$, that is 
\begin{align}
        M_{ij} = 
\begin{cases}
    A_{ij},& \text{if } \lvert A_{ij}\rvert>\text{tol or}\quad i=j\\
    0,              & \text{otherwise}
\end{cases}
\end{align}
where we set tol equal to $0.5$ and decreased it to $0.1$ as soon as the root mean square norm of the residual is close to convergence. For both the second and third options it is necessary to solve a linear system. % as given by equation \ref{eq:precond}. 
For the tridiagonal case we simply exploit a LAPACK routine which performs a Gaussian elimination with partial pivoting. Instead, in the case of the sparse $M$ matrix, the linear system is solved using the incomplete LU (iLU) decomposition\cite{gill1987maintaining} as implemented by Saunders et al\cite{ilu}.
In table \ref{tab:precond} we compare the behavior of LOBPCG and Davidson when changing the preconditioner. As expected the preconditioner based on the sparse approximation of $A$ is the one that performs best. However, such an option may become expensive and requires to store the full matrix in core or at least to have an heuristic procedure to estimate the elements. From this simple-minded experiment, we note that both methods benefit from better preconditioners, with LOBPCG exhibiting slightly more marked improvements. However, assembling such preconditioners is expensive, as it requires to build the matrix, or at least some approximation to it, in memory, which in many practical cases is far too demanding. Given the relatively small beneficial effect of going beyond a diagonal preconditioner, we believe that the latter is the optimal compromise choice.
\begin{table}[h]
\caption{Total number of iterations required to converge the lowest eigenvalue of the CASSCF orbital-rotation Hessian at the first macroiteration using the LOBPCG and Davidson solvers with three different preconditioners.}\label{tab:precond}
\begin{tabular*}{\textwidth}{@{\extracolsep\fill}ccccccc}%{@{}ccccccc@{}}
\toprule%
\multicolumn{3}{@{}c@{}}{LOBPCG} && \multicolumn{3}{@{}c@{}}{Davidson} \\\cmidrule{1-3}\cmidrule{5-7}%
diagonal & tridiagonal & iLU && diagonal & tridiagonal & iLU \\
\midrule
56 & 53 & 47 & & 52 & 54 & 50\\
\botrule
\end{tabular*}
\end{table}

\section{Conclusions and perspectives\label{sec:concl}}
In this contribution, we have described an efficient and numerically robust implementation of LOBPCG, available both in the DFTK plane-wave density functional program, and in the opensource library \textsc{diaglib}. We have discussed in detail how to avoid numerical problems and error propagation, and presented a cost-effective, yet stable strategy to orthogonalize a set of vectors using Cholesky decomposition of the overlap matrix. We have then compared the resulting implementation to Davidson's method for a selection of test-cases in quantum chemistry.
Davidson's method is the \emph{de facto} standard for solving large eigenvalue problems in quantum chemistry, and for good reasons. As many of such problems are characterized by strongly diagonally dominant matrices, Davidson's method always exhibits reliable, fast convergence. This comes, however, at a price. To be efficient, Davidson's method requires a rather large expansion subspace, which can become cumbersome for large-scale calculations, both in terms of memory requirements and computational effort in the orthogonalization step. Furthermore, the method has some difficulties dealing with non diagonally dominant matrices, as the ones encountered in CASSCF calculations.
For all these reasons, LOBPCG represents a valid alternative. Due to its low memory requirements, it can be used to treat systems for which deploying Davidson's method would be too costly. It can also be a backup method in cases where Davidson fails or, for particularly hard cases as in CASSCF, used as a default method. The implementation in \textsc{diaglib} is free and accessible, and can be used under the terms of the LGPL v2.1 license, while the one in DFTK is available under the MIT license. It is our hope that it will provide an useful tool to the developers community in quantum chemistry. 

Our numerical experiments highlight the degradation of the convergence
rate of the Davidson method as the history size is truncated. On the
other hand, LOBPCG is able to maintain a good convergence rate
(although slightly inferior to untruncated Davidson) with a subspace
of size $3N$. It would be an interesting topic of further research to
devise a method that is able to interpolate between the two, being
able to use a large history size if available, but preserving the good
behavior of LOBPCG when used with a smaller history size. This could
pave the way towards a fully adaptive method that truncates the
history size dynamically based on available information.

\section*{Acknowledgments}
This paper is dedicated to Maurizio Persico, an invaluable teacher, mentor, and colleague, and the most enthusiastic mountaineer.
A.L. wishes to thank Michael Herbst for help implementing and
stress-testing the algorithm.
F.L. T.N. and I.G. acknowledge financial support from ICSC-Centro Nazionale di Ricerca in High Performance Computing, Big Data, and Quantum Computing, funded by the European Union -- Next Generation EU -- PNRR, Missione 4 Componente 2 Investimento 1.4.
F.L. further acknowledges funding from the Italian Ministry of Research under grant 2020HTSXMA\_002 (PSI-MOVIE).

\bibliography{sn-bibliography}
\end{document}